\documentclass[a4paper, 11pt]{article}
\usepackage{mathrsfs, amsmath, amsthm}
\usepackage{subfig, graphicx, color}
\usepackage{array}
\usepackage{cases}
\makeatletter
\def\th@plain{%
  \upshape 
}
\makeatother

\makeatletter
\renewenvironment{proof}[1][\proofname]{\par
  \pushQED{\qed}%
  \normalfont \topsep6\p@\@plus6\p@\relax
  \trivlist
  \item[\hskip\labelsep
        \bfseries
    #1\@addpunct{.}]\ignorespaces
}{%
  \popQED\endtrivlist\@endpefalse
}
\makeatother

\newtheorem{theorem}{Theorem}
\numberwithin{theorem}{section}

\newtheorem{conjecture}{Conjecture}
\newtheorem*{conjecture*}{Conjecture}

\theoremstyle{definition}

\usepackage[left=25mm,top=25mm,bottom=25mm,right=25mm]{geometry}
\usepackage[T1]{fontenc}
\usepackage[varg]{txfonts}

\usepackage{enumitem}
\usepackage[square, numbers, sort&compress]{natbib}

\ifx\pdfoutput\undefined
 \usepackage[dvipdfm,%
  bookmarks=true,%
  bookmarksnumbered=true, 
  bookmarksopen=true, 
  plainpages=false,%
  pdfpagelabels,%
  colorlinks=true, 
  linkcolor=blue, 
  citecolor=blue,%
  hyperindex=true,
  urlcolor=black,
  pdfborder=001]{hyperref}
\else
 \usepackage[pdftex,%
  bookmarks=true,%
  bookmarksnumbered=true, 
  bookmarksopen=true, 
  plainpages=false,%
  pdfpagelabels,%
  colorlinks=true, 
  linkcolor=black, 
  citecolor=black,%
  anchorcolor=green,
  urlcolor= blue,
  breaklinks=true，
  hyperindex=true,
  pdfborder=001]{hyperref}
\fi

\newcounter{Hcase}
\newcounter{Hclaim}

\newcommand{\ie}{i.e.,\ }

\def\int(#1){\mathrm{int}(#1)}
\def\ext(#1){\mathrm{ext}(#1)}
\def\Int(#1){\mathrm{Int}(#1)}
\def\Ext(#1){\mathrm{Ext}(#1)}
\def\mad(#1){\mathrm{mad}(#1)}
\def\la(#1){\mathrm{la}(#1)}
\newcommand{\Lfloor}{\left\lfloor}
\newcommand{\Rfloor}{\right\rfloor}
\newcommand{\Lceil}{\left\lceil}
\newcommand{\Rceil}{\right\rceil}

\begin{document}
\title{Long properly colored cycles in edge colored complete graphs}
\author{Guanghui Wang\unskip\textsuperscript{a}, \ \  Tao Wang\unskip\textsuperscript{b, c, }\footnote{{\tt Corresponding
author: wangtao@henu.edu.cn}},
\ \ Guizhen Liu\unskip\textsuperscript{a}\\[.5em]
{\small \textsuperscript{a}\unskip School of Mathematics,}\\
{\small Shandong University, Jinan, 250100, Shandong, P. R. China}\\
{\small \textsuperscript{b}Institute of Applied Mathematics}\\
{\small Henan University, Kaifeng, 475004, P. R. China}\\
{\small \textsuperscript{c}College of Mathematics and Information Science}\\
{\small Henan University, Kaifeng, 475004, P. R. China}}
\date{}
\maketitle

\begin{abstract}
Let $K_{n}^{c}$ denote a complete graph on $n$ vertices whose edges are
colored in an arbitrary way. Let $\Delta^{\mathrm{mon}}
(K_{n}^{c})$ denote the maximum number of edges of the same color incident with a vertex of $K_{n}^{c}$. A properly colored cycle (path) in $K_{n}^{c}$ is a
cycle (path) in which adjacent edges have distinct colors. B. Bollob\'{a}s and P. Erd\"{o}s (1976) proposed the following conjecture: if $\Delta^{\mathrm{mon}}
(K_{n}^{c})<\lfloor \frac{n}{2} \rfloor$, then $K_{n}^{c}$ contains a properly colored
Hamiltonian cycle. Li, Wang and Zhou proved that if $\Delta^{\mathrm{mon}}
(K_{n}^{c})< \lfloor \frac{n}{2} \rfloor$, then $K_{n}^{c}$ contains a properly colored cycle of length at least $\lceil \frac{n+2}{3}\rceil+1$. In this paper, we improve the bound to $\lceil \frac{n}{2}\rceil + 2$.
\end{abstract}

\section{Introduction}
All graphs considered here are finite, undirected, and simple. Let $G$ be a graph with vertex set $V$ and edge set $E$.

An {\em edge coloring} of a graph is an assignment of ``colors'' to the edges of the graph. An {\em edge colored graph} is a graph with an edge coloring. A cycle (path) in an edge colored graph is {\em properly colored} if no two adjacent edges in it have the same color.

Grossman and H\"{a}ggkvist \cite{Grossman1983} gave a sufficient condition on the existence
of a properly colored cycles in edge colored graphs with two colors, and Yeo \cite{Yeo1997} extended the
result to edge colored graph with any number of colors. Given an edge colored graph $G$, let $\deg^{c}(v)$, named the color
degree of a vertex $v$, be defined as the maximum number of edges
incident to $v$ that have distinct colors. The minimum color degree $\delta^{c}(G)$ is the minimum $\deg^{c}(v)$ over all
vertices of $G$. In \cite{Li2007a},
some minimum color degree conditions for the existence of properly colored cycles
are obtained. In particular, they proved that if $G$ is an edge colored graph of order $n$ satisfying $\delta^{c}(G)\geq d\geq 2$,
 then either $G$ has a properly colored
path of length at least $2d$, or $G$ has a properly colored cycle of
length at least $\lceil\frac{2d}{3}\rceil+1$. In \cite{Lo2010}, Lo improved the bound $\lceil\frac{2d}{3}\rceil+1$ to the best possible value $d+1$.

Let $K_{n}^{c}$ be an edge colored complete graph on $n$ vertices, and let $c(u, v)$ be the color assigned to edge $uv$. In \cite{Barr1998}, Barr gave a simple sufficient condition for the existence of a
properly colored Hamiltonian path in $K_n^c$: $K_n^c$ has no monochromatic triangles. Bang-Jensen, Gutin and Yeo \cite{Bang-Jensen1998} proved that if $K_{n}^{c}$ contains a properly colored $2$-factor, then it has a properly colored Hamiltonian path. In \cite{Feng2006}, 
 Feng et al. proved  that $K^c_n$ has a properly colored Hamilton path iff $K^c_n$ has a properly colored almost $2$-factor (an almost 2-factor is a spanning subgraph consisting of disjoint cycles and a path).

Bollob\'{a}s and Erd\"{o}s \cite{Bollobas1976} proved that if $n \geq 3$ and $\delta^{c}(K_{n}^{c}) \geq \frac{7n}{8}$, then there exists a properly colored Hamiltonian cycle. They also proposed a question: whether the bound $\frac{7n}{8}$ can be improved to
$\frac{n+5}{3}$. Fujita and Magnant \cite{Fujita2011} constructed an edge colored complete graph $K_{2m}$ with $\delta^{c}(K_{2m}) = m$, which has no properly colored Hamiltonian cycles. So $\delta^{c}(K_{n}) = \lfloor\frac{n}{2}\rfloor$ is not enough. Let $\Delta^{\mathrm{mon}}(K_{n}^{c})$ be the maximum number of edges of the same color incident with a vertex of $K_{n}^{c}$. In the same paper, Bollob\'{a}s and Erd\"{o}s proposed the following conjecture.

\begin{conjecture}[Bollob\'{a}s and Erd\"{o}s \cite{Bollobas1976}]%
If $\Delta^{\mathrm{mon}}(K_{n}^{c}) < \lfloor \frac{n}{2} \rfloor$, then $K_{n}^{c}$ has a properly colored Hamiltonian cycle.
\end{conjecture}

Bollob\'{a}s and Erd\"{o}s proved that if $\Delta^{\mathrm{mon}}(K_{n}^{c}) \leq \frac{n}{69}$ then $K_{n}^{c}$ contains a properly colored Hamiltonian cycle. This result was improved by Chen and Daykin \cite{Chen1976} and Shearer \cite{Shearer1979}.

As far as we know, the best asymptotic estimate was obtained by Alon and Gutin using the probabilistic method.
\begin{theorem}[Alon and Gutin \cite{Alon1997b}]%
For every positive real number $\epsilon$, there exists $n_{0} = n_{0}(\epsilon)$ so that for every $n > n_{0}$, if $K_{n}^{c}$ satisfies
\begin{equation}%
\Delta^{\mathrm{mon}}(K_{n}^{c}) \leq (1-\frac{1}{\sqrt{2}} - \epsilon)n,
\end{equation}
then $K_{n}^{c}$ contains a properly colored Hamiltonian cycle.
\end{theorem}

Li, Wang and Zhou \cite{Li2007a} studied long properly colored cycles in edge colored complete graphs and proved that if $\Delta^{\mathrm{mon}}(K_{n}^{c}) < \lfloor \frac{n}{2} \rfloor$, then $K_{n}^{c}$ contains a properly colored cycle of length at least $\lceil \frac{n+2}{3}\rceil + 1$. For more details concerning properly colored cycles and paths,  we refer the reader to \cite{Bang-Jensen1997, Kano2008a, MR2721735, MR2519172}. In this paper, we improve the bound on the length of the properly colored cycles
and prove the following theorem.
\begin{theorem}\label{MainResult}%
If $\Delta^{\mathrm{mon}}(K_{n}^{c}) < \lfloor \frac{n}{2} \rfloor$, then $K_{n}^{c}$ contains a properly colored cycle of length at least $\lceil \frac{n}{2} \rceil + 2$.
\end{theorem}

The main idea is the rotation--extension technique of P\'{o}sa \cite{Posa1976}, which was used on general edge colored graphs in \cite{Lo2010}.

When $a$ and $b$ are integers, the notation $[a, b]$ is used to indicate the interval of all integers between $a$ and $b$, including both. In particular, $[a, a-1]$ is an empty set. A properly colored path $P$ of length $\ell$ is viewed as an $\ell$-tuple $(v_{0}, v_{1}, \dots, v_{\ell})$, which is different from $(v_{\ell}, v_{\ell-1}, \dots, v_{0})$. For such a properly colored  path $P = (v_{1}, v_{2}, \dots, v_{\ell})$, denote $\{v \mid c(v_{1}, v_{2}) \neq c(v_{1}, v)\}$ by $N^{c}(v_{1}, P)$ and $\{v \mid c(v_{\ell}, v_{\ell-1}) \neq c(v, v_{\ell})\}$ by $N^{c}(v_{\ell}, P)$.
\section{The Proof of Theorem 1.2}

If $3 \leq n \leq 5$, then $K_{n}^{c}$ is a properly colored complete graph, and hence each Hamiltonian cycle is properly colored. Our conclusion holds clearly. So we may assume that $n \geq 6$. 

We will prove Theorem \ref{MainResult} by contradiction. Suppose that each properly colored cycle is of length at most $\lceil \frac{n}{2} \rceil + 1$. For simplicity, let the vertices of $K_n^c$ be labeled with integers from $1$ to $n$. Let $P_{0}$ be a longest properly colored path. Without loss of generality, assume that $P_0=(1, 2, \dots, \ell)$.

Put $$X = \{x_{1}, x_{2}, \dots, x_{p}\} = \{x_{i} \mid c(1, x_{i}) \neq c(1, 2)\}=N^{c}(1, P_{0})$$ and $$Y = \{y_{1}, y_{2}, \dots, y_{q}\} = \{y_{i} \mid c(y_{i}, \ell) \neq c(\ell-1, \ell)\}=N^{c}(\ell, P_{0}).$$

Moreover, $x_{1}, x_{2}, \dots, x_{p}$ and $y_{1}, y_{2}, \dots, y_{q}$ are increasing sequences. Since $\Delta^{\mathrm{mon}}(K_{n}^{c}) < \lfloor \frac{n}{2} \rfloor$, we have $\min\{|X|, |Y|\} \geq \lceil \frac{n}{2} \rceil$; consequently, $x_{1} < x_{p}$. Note that $X, Y\subset V(P_0)$ since $P_0$ is a longest properly colored path. Thus $\ell \geq \lceil \frac{n}{2} \rceil + 2$. Note that either $\ell\notin X$ or $1\notin Y$; otherwise $(1,2,\dots, \ell,1)$ is a properly colored cycle of length at least $\lceil \frac{n}{2} \rceil +2$, which is a contradiction. Hence $\ell \geq \lceil \frac{n}{2} \rceil + 3$.

Let $y_{s} \in Y$ be the maximum such that $c(\ell, y_{i}) = c(y_{i}, y_{i}+1)$ for all $y_{i} \in [y_{1}, y_{s}] \cap Y$. Note that $y_{s}$ is well defined, since $c(\ell, y_{1}) = c(y_{1}, y_{1}+1)$; otherwise, $(y_{1}, y_{1}+1, \dots, \ell, y_{1})$ is a properly colored cycle of length at least $\lceil \frac{n}{2} \rceil + 2$. Clearly, $y_{s} \leq \ell-2$.

Note that $c(1, x_{p}) = c(x_{p}-1, x_{p})$ or else $(1, 2, \dots, x_{p}, 1)$ is a properly colored cycle of length at least $\lceil \frac{n}{2} \rceil + 2$, a contradiction.

{\noindent\bf Claim 1.} $y_{s} \leq x_{p}-3$.

\begin{proof}[Proof of Claim 1]
If $x_{p} = \ell$, then $y_{s} \leq \ell-3$; otherwise, $(1, 2, \dots, \ell-2, \ell, 1)$ is a properly colored cycle of length $\ell-1\geq \lceil \frac{n}{2} \rceil + 2$, a contradiction. So we may assume that $x_{p} \leq \ell-1$. Let $y_{i} \in Y$ be the maximum such that $y_{i} < x_{p}$. If $c(y_{i}, y_{i}+1) = c(y_{i}, \ell)$, then $(1, \dots, y_{i}, \ell, \ell-1, \dots, x_{p}, 1)$ is a properly colored cycle containing $Y \cup \{\ell, \ell-1\}$, a contradiction. Hence, $c(y_{i}, y_{i}+1) \neq c(y_{i}, \ell)$, by the definition of $y_{s}$, we have $y_{s} \leq x_{p}-2$. If $y_{s} = x_{p} - 2$, then $(1, \dots, x_{p}-2, \ell, \ell - 1, \dots, x_{p}, 1)$ is a properly colored cycle of length $\ell-1 \geq \lceil \frac{n}{2} \rceil + 2$, which is a contradiction. Therefore, $y_{s} \leq x_{p}-3$.
\end{proof}

{\noindent\bf Claim 2.} There exist $u, w \in X$ such that $1 \leq y_{1} \leq y_{s} < u \leq \lceil \frac{n}{2} \rceil + 1$, $u < w$ and the following proposition holds:
\begin{enumerate}[label= (\alph*)]
    \item $c(y_{i}, \ell) = c(y_{i}, y_{i}+1)$ for all $y_{i} \in [y_{1}, y_{s}] \cap Y$;
    \item $c(1, x_{i}) = c(x_{i}, x_{i}+1)$ for all $x_{i} \in [y_{1}+1, u] \cap X$;
    \item $c(1, w) \neq c(w, w+1)$ or $w = \ell$;
    \item if $x_{i} \in X$ and $x_{i} < w$, then $x_{i} \leq u$.
\end{enumerate}
\begin{proof}[Proof of Claim 2]
(a) By the definition of $y_{s}$, the assertion (a) holds.

(b) First, we show that $c(1, x_{i}) = c(x_{i}, x_{i}+1)$ for all $x_{i} \in [y_{1}+1, y_{s}+1] \cap X$. Otherwise, if there exists $x_{i} \in [y_{1}+1, y_{s}+1] \cap X$ such that $c(1, x_{i}) \neq c(x_{i}, x_{i}+1)$, let $y_{j}$ be the maximum such that $y_{j} \in [y_{1}, y_{s}] \cap Y$ and $y_{j} < x_{i}$, then $(1, 2, \dots, y_{j}, \ell, \ell-1, \dots, x_{i}, 1)$ is a properly colored cycle containing $Y \cup \{\ell-1, \ell\}$, a contradiction.

Next, we show that there exists $u$ such that $c(1, x_{i}) = c(x_{i}, x_{i}+1)$ for all $x_{i} \in [y_{s}+1, u] \cap X$. Let $x \in X$ be the minimum such that $x > y_{s}$. Note that $x$ must exist since $x_{p} > y_{s}$ by Claim 1. If $x = \ell$, then $y_{s} \geq y_{1} \geq 2$; hence $(1, \dots, y_{s}, \ell, 1)$ is a properly colored cycle containing $X \cup \{1, 2\}$, a contradiction, so $x < \ell$. Suppose that $c(1, x) \neq c(x, x+1)$. If $y_{s} = 1$, then $(1, \ell, \ell-1, \dots, x, 1)$ is a properly colored cycle containing $X \cup \{1, \ell\}$, a contradiction; if $y_{s} > 1$, then $(1, \dots, y_{s}, \ell, \ell-1, \dots, x, 1)$ is a properly colored cycle containing $X \cup \{1, 2\}$, a contradiction. So, we have $x < \ell$ and $c(1, x) = c(x, x+1)$.

Let $u$ be the maximum such that $c(1, x_{i}) = c(x_{i}, x_{i}+1)$ for all $x_{i} \in [y_{s}+1, u] \cap X$ and $y_{s} < u < \ell$. By the above argument, $u$ is well defined. In particular, $c(1, u) = c(u, u+1)$, it follows that $(1, 2, \dots, u, 1)$ is a properly colored cycle, and hence $u \leq \lceil \frac{n}{2} \rceil + 1$.

Let $w \in X$ be the minimum such that $w>u$. Since $c(1, x_{p}) = c(x_{p} - 1, x_{p}) \neq c(x_{p}, x_{p} +1)$, it follows that the vertex $w$ must exist, and then (c) and (d) hold.
\end{proof}

Let $S(P_{0}) = Y \cap [y_{1}, y_{s}]$. Note that $y_{1}, y_{s}, u, w, S$ are viewed as functions of longest properly colored paths.

By Claims 2~(a) and 2~(c), $C_{0}=(1, \dots, y_{s}, \ell, \ell-1, \dots, w, 1)$ is a properly colored cycle; see Fig.~1.

\begin{figure}[!htb]
 \centering
\includegraphics{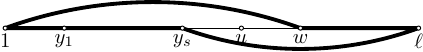}\\
\caption{Cycle $C_{0}=(1, \dots, y_{s}, \ell, \ell-1, \dots, w, 1)$}
\label{}
\end{figure}

Without loss of generality, assume that $P_{0}$ is a longest properly colored path satisfying that $|S(P_{0})|$ is maximum over all the longest properly colored paths.

Clearly, $P^{*} = (y_{1}+1, y_{1}+2, \dots, \ell, y_{1}, \dots, 1)$ is a longest properly colored path. Note that $$N^{c}(1, P^{*})=X.$$
Since $c(y_{1}+1, y_{1}+2) \neq c(y_{1}+1, y_{1})$, it follows that $y_{1} \in N^{c}(y_{1}+1, P^{*})$. Hence $P^{*}$ is a longest properly colored path.

{\noindent\bf Claim 3.} $|X| = \lceil \frac{n}{2} \rceil$ and $|C_{0}| = \lceil \frac{n}{2} \rceil + 1$. Moreover, $S(P_{0}) = [y_{1}, y_{s}]$, \ie it is an integer interval, and
\begin{numcases}{X=}
[3, y_{1}] \cup [t, u] \cup [w, \ell]     & if $y_{1} > 1$;\\{}
[t, u] \cup [w, \ell-1]        & if $y_{1} = 1$,
\end{numcases}
where $|[t, u]| = |[y_{1}, y_{s}]|$, $t \geq \max\{y_{1}+1, 3\}$ and $S(P^{*}) = [t, u]$.

\begin{proof}[Proof of Claim 3]
By Claim 2~(b), $[y_{1}+1, u] \cap X \subseteq S(P^{*})$.
Since $c(1, u) = c(u, u+1)$, it yields that $u \in S(P^{*})$. If $w < \ell$, then $c(1, w) \neq c(w, w+1)$, and then $w \notin S(P^{*})$; if $w=\ell$, then $c(1, \ell) \neq c(\ell, y_{1})$, and then $w \notin S(P^{*})$. Therefore, $S(P^{*}) = [y_{1}+1, u] \cap X$ and $|[y_{1}, y_{s}] \cap Y| \geq |[y_{1}+1, u] \cap X|$ by the maximality of $|S(P_{0})|$.

Now
\begin{align*}%
\Lceil \frac{n}{2} \Rceil +1
&\geq |C_{0}|\\
& =|[1, y_{s}] \cap X| + |[w, \ell] \cap X| + |[1, y_{s}]\setminus X| + |[w, \ell]\setminus X|\\
& =|X| - |[y_{s}+1, u] \cap X| + |[1, y_{s}]\setminus X| + |[w, \ell]\setminus X|\\
& =|X| - |[y_{1}+1, u] \cap X| + |[y_{1}+1, y_{s}] \cap X| + |[1, y_{s}] \setminus X| + |[w, \ell] \setminus X|\\
&\geq |X| - |[y_{1}, y_{s}] \cap Y| + |[y_{1}+1, y_{s}]| + |[1, y_{1}]\setminus X| + |[w, \ell] \setminus X|\\
&= |X| - |[y_{1}, y_{s}] \cap Y| + |[y_{1}, y_{s}]| + |[2, y_{1}]\setminus X| + |[w, \ell] \setminus X|\\
&\geq |X| + |[2, y_{1}]\setminus X| + |[w, \ell] \setminus X|\\
&\geq |X|+1\\
&\geq \Lceil \frac{n}{2} \Rceil +1.
\end{align*}

Therefore, all the inequalities become equalities. Then $|C_{0}| = \lceil \frac{n}{2} \rceil +1$, $|X| = \lceil \frac{n}{2} \rceil$, $|S(P_{0})| = |S(P^{*})|$, $S(P_{0}) = [y_{1}, y_{s}]$ and
\begin{equation}\label{EQ1}%
|[2, y_{1}]\setminus X| + |[w, \ell] \setminus X|= 1.
\end{equation}
\eqref{EQ1} implies that
\begin{numcases}{}
[3, y_{1}] \subseteq X\mbox{ and }[w, \ell] \subseteq X & if $y_{1} > 1$;\\{}
[w, \ell] \cap X = [w, \ell-1] & if $y_{1} = 1$.
\end{numcases}

Since $|S(P_{0})| = |S(P^{*})|$, by the above argument, $S(P^{*})$ must be an integer interval, so $S(P^{*}) = [y_{1}+1, u] \cap X= [t, u]$, clearly $t \geq \max\{y_{1}+1, 3\}$.
\end{proof}

Note that $t$ is also viewed as a function of the longest properly colored path.

{\noindent\bf Claim 4.}
\begin{enumerate}[label = (\alph*)]
 \item For any $y \in Y \cap [y_{s}+1, w-1]$, we have $c(y, y-1) = c(y, \ell)$; $|S(P_{0})| \geq |Y \cap [y_{s} + 1, w - 1]|$.
 \item If $w = x_{p}$, then $X = [3, \lceil \frac{n}{2} \rceil +1] \cup \{x_{p}\}$.
\end{enumerate}

\begin{proof}[Proof of Claim 4]
(a) Since $|C_{0}| = \lceil \frac{n}{2} \rceil +1$, we have $c(y, y-1) = c(y, \ell)$ for any $y \in Y \cap [y_{s}+1, w-1]$; otherwise, $(1, 2, \dots, y, \ell, \ell-1, \dots, w, 1)$ is a properly colored cycle of length greater than $|C_{0}|$, a contradiction.

Obviously, $Q=(w-1, w-2, \dots, 1, w, w+1, \dots, \ell)$ is a longest properly colored path. Note that $N^{c}(\ell, Q) = Y$. We have just proved $Y \cap [y_{s}+1, w-1] \subseteq S(Q)$, and hence $|S(P_{0})| \geq |S(Q)| \geq Y \cap [y_{s}+1, w-1]$.

(b) Since $w = x_{p}$, it follows that $u = x_{p-1}$ and $(1, 2, \dots, x_{p-1}, 1)$ is a properly colored cycle containing $\{1, 2, x_{1}, \dots, x_{p-1}\}$, so $x_{p-1} = \lceil \frac{n}{2} \rceil +1$. Hence $X = [3, \lceil \frac{n}{2} \rceil +1] \cup \{x_{p}\}$.
\end{proof}

{\noindent\bf Claim 5.} $|S(P_{0})| \geq 3$.

\begin{proof}[Proof of Claim 5]
By way of contradiction, we assume that $|S(P_{0})| \leq 2$. Without loss of generality, we may assume that $y_{1} \neq 1$, otherwise consider the properly colored path $(2, 3, \dots, \ell, 1)$ instead.

{\noindent\bf Case 1:} $w \leq \ell-1$. 

Let $P'=(\ell, \ell-1, \dots, 1)$, it is a longest properly colored path. If $w \leq \ell-1$, then $c(1, \ell) = c(\ell, \ell-1)$ and $c(1, \ell-1)= c(\ell-1, \ell-2)$, so $\ell, \ell-1 \in S(P')$. Hence, $|S(P_{0})| =|S(P')| = 2$ and $c(1, \ell-2) \neq c(\ell-2, \ell-3)$. Applying Claim 3 with $P'$, we have $c(\ell, \ell-1) \neq c(\ell, 2)$. Notice that $c(\ell, 2) = c(2, 3)$, for otherwise $(2, 3, \dots, \ell, 2)$ is a properly colored cycle of length at least $\lceil \frac{n}{2} \rceil + 2$, a contradiction. So we have that $y_{1} = 2$ and $y_{s} = 3$.  

{\noindent\bf Case 1.1:} $w \leq \ell-2$. Then $C_1=(1, 2, \dots, \ell-2, 1)$ is a properly colored cycle of length at most $\lceil \frac{n}{2} \rceil + 1$, \ie $\ell-2 \leq \lceil \frac{n}{2} \rceil +1$, or $\ell \leq \lceil \frac{n}{2} \rceil +3$. On the other hand, $\ell \geq \lceil \frac{n}{2} \rceil +3 $, and then $\ell= \lceil \frac{n}{2} \rceil +3$. Moreover, $Y=[2,\ell-2]$. Since the properly colored cycle $C_0=(1,2,3,\ell,\ell-1,\dots,w,1)$ has length $\lceil \frac{n}{2} \rceil +1$ by Claim 3, we conclude that $w=6$. Note that $|X| = \lceil \frac{n}{2} \rceil = |[t, u]| + |[w, \ell ]| \geq 5$ implying $n \geq 9$. 

If $\ell = n$, then $n-2 = \ell -2 = \lceil \frac{n}{2} \rceil +1$ and $n \in \{6, 7\}$, a contradiction. Then we may assume that there exists a vertex $z$ such that it is not on $P_{0}$. The properly colored path $(1, 2, 3, \ell, \ell-1, \dots, 5, 4)$ is longest, it follows that $c(4, 5) = c(4, z)$; otherwise, $(1, 2, 3, \ell, \ell-1, \dots, 5, 4, z)$ is a longer properly colored path. By Claim~4(a), we have $c(\ell, 5) = c(5, 4) \neq c(5, 6)$. Therefore, $(\ell-1, \ell, 5, 6, \dots, \ell-2, 1, 2, 3, 4, z)$ is a longer properly colored path than $P_{0}$, which is a contradiction.

{\noindent\bf Case 1.2:} $w = \ell-1$.  Then $|X| = \lceil \frac{n}{2} \rceil = |[t, u]| + |[w, \ell]| = 4$, \ie $n \in \{7, 8\}$. Since $\Lceil \frac{n}{2} \Rceil +3\leq \ell\leq n$, $\ell=7$ or $\ell=8$.

If $t = 3$, then $c(1, 3)=c(3, 4)$. Since $y_{s} = 3$, $c(3, 4)=c(\ell, 3)$. Consider the properly colored path $P'' = (3, 4, \dots, \ell, 2, 1)$. By Claim 3, we have $S(P'') = \{3, 4\}$ and $N^{c}(3, P'') = \{t', t'+1\} \cup \{\ell, 2\}$, so $c(\ell, 3) \neq c(3, 4)$, a contradiction.

If $t = 4$, then $X = \{4, 5, \ell-1, \ell\}$. Define $\phi$ be the permutation on $[1, \ell]$ such that $(\phi(1), \phi(2), \dots, \phi(\ell)) = (3, 4, \dots, \ell, 2, 1)$.
We will show the following propositions for $i\geq 0$ by induction on $i$.
\begin{enumerate}[label = (\alph*)]
 \item $P_{i} = (\phi^{i}(1), \phi^{i}(2), \dots, \phi^{i}(\ell))$ is a longest properly colored path;
 \item $S(P_{i})=\{\phi^{i}(2), \phi^{i}(3)\}$;
 \item $N^{c}(\phi^{i}(1), P_{i}) = \{\phi^{i}(j) \mid j \in [4, 5] \cup [\ell - 1, \ell]\}$.
 \item $c(\phi^{i}(1), \phi^{i}(j)) = c(\phi^{i}(j), \phi^{i}(j+1))$ for $j \in [4, 5]$.
\end{enumerate}

Clearly, the propositions (a)--(d) hold for $i = 0$. Assume that $i \geq 1$ and they are true for $i-1$. Let $j'$ denote $\phi^{i-1}(j)$. Since $2'\in S(P_{i-1})$, $c(1',2')\neq c(2',3')=c(2',\ell')$. Hence, $P_{i} = (\phi(1'),\phi(2'),\cdots,\phi(\ell'))$ is a longest properly colored path. Also, by the definition of $N^{c}(\phi^{i-1}(1), P_{i-1})$ and $|S(P_{i})| = |S(P_0)| = 2$, (b) holds. By Claim 3 and $|X| = 4$,

$$N^{c}(\phi^i(1), P_{i}) = \{\phi^{i}(j) \mid j \in [t_{i}, t_{i} + 1] \cup [\ell-1, \ell]\},$$

\noindent for some $t_{i} \in \{3, 4\}$. If $t_{i} = 3$, we can also get a contradiction as above by taking $P_{0} = P_{i}$. So $t_{i} = 4$, and thus (c) holds. By Claim 3, (d) also holds.

If $\ell= 7$, then $X = \{4, 5, 6, 7\}$ and $Y = \{2, 3, 4, 5\}$. For the path $P_{0}$, we have $c(7, 2) = c(2, 3)$ and $c(7, 3) = c(3, 4)$. Taking $i = 1$ in the proposition, consider the longest properly colored path $P_{1}$, \ie $(3, 4, 5, 6, 7, 2, 1)$. By proposition (c), $c(\phi^{1}(1), \phi^{1}(5)) = c(\phi^{1}(5), \phi^{1}(6))$, \ie $c(3, 7) = c(7, 2)$. Therefore, $c(2, 3) = c(7, 2) = c(3, 7) = c(3, 4)$, a contradiction.

If $\ell= 8$, then $X = \{4, 5, 7, 8\}$. Since $u = 5$, it follows that $c(1, 5) = c(5, 6) \neq c(4, 5)$. Consider the properly colored path $P_{5} = (4, 3, 6, 5, 8, 7, 1, 2)$. Since $8 = \phi^{5}(5)$, by proposition (c) and Claim 2~(b), $c(4, 8) = c(8, 7) \neq c(5, 8)$; moreover, $5 \in N^{c}(8, P_{0})$. Thus $c(8, 5) = c(4, 5)$ by Claim 4. Hence $c(1, 5) = c(5, 6) \neq c(4, 5) = c(5, 8)$. Now, we have a properly colored cycle $(1, 2, 3, 4, 8, 5, 1)$ of length 6, a contradiction.

{\noindent\bf Case 2:} $w=\ell$. Since $|Y \cap [y_{s} + 1, w - 1]| \leq |S(P_{0})|$ by Claim 4~(a), it follows that $|Y| = |[y_{1}, y_{s}]| + |Y \cap [y_{s} + 1, w - 1]| \leq 2|S(P_{0})|$. If $|S(P_{0})| = 1$, then $\lceil \frac{n}{2} \rceil \leq |Y| \leq 2$, and thus $n \leq 4$, a contradiction. Hence, $|S(P_{0})| = 2$, and $\lceil \frac{n}{2} \rceil \leq |Y| \leq 4$, so $n \leq 8$. By Claim 3, $X = [3, \lceil \frac{n}{2} \rceil +1] \cup \{\ell\}$, $u=\lceil \frac{n}{2} \rceil +1$, $t=\lceil \frac{n}{2} \rceil$, $y_{1} = t-1= \lceil \frac{n}{2} \rceil - 1$. So $Y \subseteq [\lceil \frac{n}{2} \rceil - 1, \ell-2]$; hence, $\lfloor \frac{n}{2} \rfloor = n-2 - (\lceil \frac{n}{2} \rceil - 1) + 1 \geq \ell-2 - (\lceil \frac{n}{2} \rceil - 1) + 1 \geq |Y| \geq \lceil \frac{n}{2} \rceil$, and then $n=\ell$ is even, $Y = [ \frac{n}{2}  - 1, n-2]$ and $|Y| = \frac{n}{2}$. Now, we know that $n \in \{6, 8\}$.

If $n = 6$, then $X = \{3, 4, 6\}$, $Y = \{2, 3, 4\}$, $y_{1} = 2$, $y_{s} = t = 3$ and $u = 4$. By Claim 2~(a) and 2~(b), $c(1, 3) = c(3, 4) = c(3, 6)$, which contradicts the fact that $\Delta^{\mathrm{mon}}(K_{6}^{c}) < 3$.

If $n = 8$, then $X = \{3, 4, 5, 8\}$, $Y = \{3, 4, 5, 6\}$, $y_{1} = 3$, $t = 4$ and $u = 5$. Since $4 \in [y_{1}+ 1, u] \cap X$, $c(1, 4) = c(4, 5) \neq c(3, 4)$ by Claim 2~(b); by Claim 2~(a), $c(8, 3) = c(3, 4) \neq c(2, 3)$ and $c(8, 4) = c(4, 5) \neq c(3, 4)$. Hence, $c(4, 1) = c(4, 5) = c(4, 8)$. Since $\Delta^{\mathrm{mon}}(K_{8}^{c}) < 4$, it follows that $c(4, 7) \neq c(4, 1)$. If $c(4, 7) \neq c(7, 8)$, then $(1, 2, 3, 8, 7, 4, 1)$ is a properly colored cycle of length 6, a contradiction. So $c(4, 7) = c(7, 8) \neq c(7, 6)$. By Claim 4~(a), we have $c(8, 5) = c(4, 5) \neq c(5, 6)$, so $c(8, 5) \neq c(8, 3)$. Now we have a properly colored Hamiltonian cycle $(1, 2, 3, 8, 5, 6, 7, 4, 1)$, a contradiction.
\end{proof}

{\noindent\bf Claim 6.} Without loss of generality, we may assume that $t \geq y_{1} + 3$.

\begin{proof}[Proof of the Claim 6]%
Suppose to the contrary that $t \leq y_{1} + 2$. Then $t=y_1+1$ or $t=y_1+2$. Since $|S(P_0)|\geq 3$, we have $t-1\in S(P_0)$. Without loss of generality, we may assume that $y_{1} = 1$; otherwise, consider the longest properly colored path $(t, t+1, \dots, \ell, t-1, t-2, \dots, 1)$ instead. By Claim 3, $t = 3$. Moreover, $X = [3, y_{s}+2] \cup [w, \ell - 1]$.

Let $\phi$ be a permutation on $[1, \ell]$ such that $(\phi(1), \phi(2), \dots, \phi(\ell-2), \phi(\ell-1), \phi(\ell)) = (3, 4, \dots, \ell, 2, 1)$. We will show the following statements for $i \geq 0$ by induction.
\begin{enumerate}[label = (\alph*)]
 \item $P_{i} = (\phi^{i}(1), \phi^{i}(2), \dots, \phi^{i}(\ell))$ is a longest properly colored path;
 \item $S(P_{i}) = \{\phi^{i}(j) \mid j \in [1, y_{s}]\}$;
 \item $N^{c}(\phi^{i}(1), P_{i}) = \{\phi^{i}(j) \mid j \in [3, y_{s} + 2] \cup [w, \ell - 1]\}$;
 \item for all $j \in [3, y_{s}+2]$, we have $c(\phi^{i}(1), \phi^{i}(j)) = c(\phi^{i}(j), \phi^{i}(j+1))$;
 \item for all $j \in [4, \lceil \frac{n}{2} \rceil + 2]$, we have $c(\phi^{i}(2), \phi^{i}(j)) = c(\phi^{i}(j), \phi^{i}(j+1))$.
\end{enumerate}
First, the statements are true for $i = 0$. Assume that the statements are true for $i - 1$, where $i \geq 1$. Since $\phi^{i-1}(2) \in S(P_{i-1})$, it follows that $P_{i} = (\phi^{i-1}(3), \phi^{i-1}(4), \dots, \phi^{i-1}(\ell), \phi^{i-1}(2), \phi^{i-1}(1)) = (\phi^{i}(1), \phi^{i}(2), \dots, \phi^{i}(\ell))$ is a longest properly colored path. Since $|S(P_{i-1})| = |S(P_{0})| = |[1, y_{s}]|$, we have $t(P_{i-1}) = \phi^{i-1}(3)$; by Claim 3, it yields that $S(P_{i}) = [t(P_{i-1}), u(P_{i-1})] = [\phi^{i}(1), \phi^{i}(y_{s})]$, and thus $|S(P_{i})|=|S(P_0)|$. Hence, by the assumption, we have $N^{c}(\phi^{i}(1), P_{i}) = \{\phi^{i}(j) \mid j \in [3, y_{s} + 2] \cup [w, \ell - 1]\}$. By Claim 2~(b), for all $j \in [3, y_{s} + 2]$, $c(\phi^{i}(1), \phi^{i}(j)) = c(\phi^{i}(j), \phi^{i}(j+1))$.

Since $S(P_{i}) = \{\phi^{i}(j) \mid j \in [1, y_{s}]\}$, it follows that $P_{i}' = (\phi^{i}(2), \phi^{i}(3), \dots, \phi^{i}(\ell), \phi^{i}(1))$ is a longest properly colored path. By Claim 3, $S(P_{i}') = \{\phi^{i}(j) \mid j \in [3, y_{s} + 2]\}$. By Claim 2~(a), $c(\phi^{i}(2), \phi^{i}(3)) = c(\phi^{i}(2), \phi^{i}(\ell))$, so $\phi^{i}(\ell) \notin N^{c}(\phi^{i}(2), P_{i}')$; by Claim 3, $N^{c}(\phi^{i}(2), P_{i}') = \{\phi^{i}(j) \mid j \in [4, \lceil \frac{n}{2}\rceil + 2] \cup \{1\}\}$. By Claim 2~(b), $c(\phi^{i}(2), \phi^{i}(j)) = c(\phi^{i}(j), \phi^{i}(j+1))$ for all $j \in [4, \lceil \frac{n}{2} \rceil + 2]$. This completes the proof of the statements.

Consider the properly colored path $P_{0}$. Since $|S(P_{0})| \geq 3$, it follows that $y_{s} \geq 3$. Thus $5 \in [3, y_{s} + 2]$. By statement (d), $c(1, 5) = c(5, 6)$.

If $\ell$ is odd, take $i = \frac{\ell+1}{2}$. It is easy to check that $P_{\frac{\ell+1}{2}} = (1, 4, 3, 6, 5, \dots)$. By statement (d), $c(1, 5) \neq c(6, 5)$, a contradiction.

If $\ell$ is even, take $i = \frac{\ell}{2}$. It is easy to check that $P_{\frac{\ell}{2}} = (2, 1, 4, 3, 6, 5, \dots)$. By statement (e), if $\lceil \frac{n}{2} \rceil + 2 \geq 6$, then $c(1, 5) = c(5, 8) \neq c(5, 6)$, a contradiction. So $n = 6$. Since $3 \in [1, y_{s}]$, it yields that $c(3, 4) = c(3, 6)$ by Claim 2~(a); by Claim 2~(b), $c(1, 3) = c(3, 4)$. Therefore, $c(1, 3) = c(3, 4) = c(3, 6)$, which contradicts the fact that $\Delta^{\mathrm{mon}}(K_{6}^{c}) < \lfloor \frac{6}{2} \rfloor$.
\end{proof}

{\noindent\bf Claim 7.} $c(y_{1} + 1, y_{1} + 3) \notin \{c(y_{1} + 1, y_{1} + 2), c(y_{1} + 3, y_{1} + 4)\}$.
\begin{proof}[Proof of Claim 7]%
By Claim 3, $S(P^{*}) = [t, u]$ and $y_{1} + 3 \in N^{c}(y_{1}+1, P^{*})$, so $c(y_{1} + 1, y_{1} + 2) \neq c(y_{1} + 1, y_{1} + 3)$. For the longest properly colored path $Q^{*} = (y_{1} + 3, y_{1} + 4, \dots, \ell, y_{1} + 2, y_{1} + 1, \dots, 1) = (\psi(1), \psi(2), \dots, \psi(\ell))$, we know that $S(Q^{*}) = [t, u]$ and $N^{c}(y_{1} + 3, Q^{*}) = \{\psi(j) \mid j \in A\}$, where
\begin{numcases}{A=}
[3, y_{1}''] \cup [t'', u''] \cup [w'', \ell]     & if $t > y_{1} + 3$;\\{}
[t'', u''] \cup [w'', \ell-1]        & if $t = y_{1} + 3$,
\end{numcases}
for some $2 \leq y_{1}'' < t'' < u'' < w'' \leq \ell$. Since $y_{1} + 2 \in S(P_{0})$, $c(y_{1} + 2, y_{1} + 3) = c(y_{1} + 2, \ell)$. Thus $y_{1} + 2 \in \{\psi(j) \mid j \in [w'', \ell-1]\}$, so $[2, y_{1} + 2] \subseteq \{\psi(j) \mid j \in [w'', \ell-1]\}$. In particular, $y_{1} + 1 \in N^{c}(y_{1} + 3, Q^{*})$, so $c(y_{1} + 1, y_{1} +3) \neq c(y_{1} + 3, y_{1} + 4)$, which completes the proof of Claim 7.
\end{proof}

\begin{figure}[!htb]
 \centering
 \includegraphics{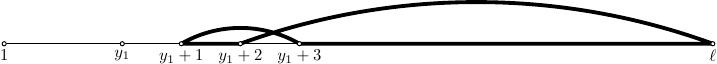}
 \caption{Cycle $(y_{1} + 1, y_{1} +2, \ell, \ell - 1, \dots, y_{1} + 3, y_{1} +1)$}
 \label{}
\end{figure}

By Claim 7, note that $(y_{1} + 1, y_{1} +2, \ell, \ell - 1, \dots, y_{1} + 3, y_{1} +1)$ is a properly colored cycle (see Fig.~2) containing $Y \cup \{\ell-1, \ell\} \setminus \{y_{1}\}$, and then its length is exactly $\lceil \frac{n}{2} \rceil + 1$. Moreover, $Y= [\ell - \lceil \frac{n}{2} \rceil -1, \ell-2]$.

If $\ell = n$, then $Y = [\lfloor \frac{n}{2} \rfloor - 1, n-2]$ and $y_{1} = \lfloor \frac{n}{2} \rfloor - 1$. By Claim 3 and Claim 6, $u \geq y_{s} + 3$. Now

\begin{equation}\label{EQ2}%
\Lceil \frac{n}{2} \Rceil +1 \geq u \geq y_{s} + 3 = y_{1} + |S(P_{0})| + 2 \geq \Lfloor \frac{n}{2} \Rfloor + 4,
\end{equation}
a contradiction.

Then we may assume that $\ell < n$, \ie there exists a vertex $z$ which is not on $P_{0}$. We know that $c(\ell-1, \ell) = c(\ell, z)$; otherwise $(1, 2, \dots, \ell, z)$ is a properly colored path which is longer than $P_{0}$. We have known that $(1, 2, \dots, y_{1}+1, \ell, \ell-1, \dots, y_{1} + 2)$ is a longest properly colored path; it yields that $c(y_{1} + 2, y_{1} + 3) = c(y_{1} + 2, z)$ for the same reason. Since $y_{1} + 2 \in S(P_{0})$, we have $c(y_{1}+2, y_{1}+3) = c(y_{1}+2, \ell) \neq c(\ell-1, \ell) = c(\ell, z)$, and hence $c(y_{1}+2, z) \neq c(\ell, z)$. Therefore, $(1, \dots, y_{1}, \ell, z, y_{1}+2, y_{1}+1, y_{1}+3, \dots, \ell-1)$ is a properly colored path which is longer than $P_{0}$, a contradiction.

\vskip 3mm \vspace{0.3cm} \noindent{\bf Acknowledgments.} This research was supported by NSFC Grants (61373027, 11101243, 11101125), NSF of Shandong Province (ZR2012FM023) and the Scientific Research Foundation for the Excellent Middle-Aged and Young Scientists of Shandong Province of China (BS2012SF016). In addition, the authors would like to thank the anonymous reviewers for their valuable comments and assistance on earlier drafts.


\begin{thebibliography}{10}

\bibitem{Alon1997b}
N.~Alon and G.~Gutin, Properly colored {Hamilton} cycles in edge-colored
  complete graphs, Random Structures and Algorithms 11 (1997)~(2) 179--186.

\bibitem{Bang-Jensen1997}
J.~Bang-Jensen and G.~Gutin, Alternating cycles and paths in edge-coloured
  multigraphs: A survey, Discrete Math. 165-166 (1997) 39--60.

\bibitem{Bang-Jensen1998}
J.~Bang-Jensen, G.~Gutin and A.~Yeo, Properly coloured {Hamiltonian} paths in
  edge-coloured complete graphs, Discrete Appl. Math. 82 (1998)~(1-3) 247--250.

\bibitem{Barr1998}
O.~Barr, Properly coloured {Hamiltonian} paths in edge-coloured complete graphs
  without monochromatic triangles, Ars Combin. 50 (1998) 316--318.

\bibitem{Bollobas1976}
B.~Bollob\'{a}s and P.~Erd\"{o}s, Alternating hamiltonian cycles, Israel J.
  Math. 23 (1976)~(2) 126--131.

\bibitem{Chen1976}
C.~C. Chen and D.~E. Daykin, Graphs with {Hamiltonian} cycles having adjacent
  lines different colors, J. Combin. Theory Ser. B 21 (1976)~(2) 135--139.

\bibitem{Feng2006}
J.~Feng, H.-E. Giesen, Y.~Guo, G.~Gutin, T.~Jensen and A.~Rafiey,
  Characterization of edge-colored complete graphs with properly colored
  {Hamilton} paths, J. Graph Theory 53 (2006)~(4) 333--346.

\bibitem{Fujita2011}
S.~Fujita and C.~Magnant, Properly colored paths and cycles, Discrete Appl.
  Math. 159 (2011)~(14) 1391--1397.

\bibitem{Grossman1983}
J.~W. Grossman and R.~H\"{a}ggkvist, Alternating cycles in edge-partitioned
  graphs, J. Combin. Theory Ser. B 34 (1983)~(1) 77--81.

\bibitem{MR2721735}
G.~Gutin and E.~J. Kim, Properly coloured cycles and paths: results and open
  problems, Lecture Notes in Comput. Sci. 5420 (2009) 200--208.

\bibitem{Kano2008a}
M.~Kano and X.~Li, Monochromatic and heterochromatic subgraphs in edge-colored
  graphs - a survey, Graphs Combin. 24 (2008)~(4) 237--263.

\bibitem{Li2007a}
H.~Li, G.~Wang and S.~Zhou, Long alternating cycles in edge-colored complete
  graphs, Lecture Notes in Comput. Sci. 4613 (2007)
  305--309.

\bibitem{Lo2010}
A.~Lo, A {Dirac} type condition for properly coloured paths and cycles, J.
  Graph Theory (2013) http://dx.doi.org/10.1002/jgt.21751.

\bibitem{Posa1976}
L.~P\'{o}sa, Hamiltonian circuits in random graphs, Discrete Math. 14
  (1976)~(4) 359--364.

\bibitem{Shearer1979}
J.~Shearer, A property of the colored complete graph, Discrete Math. 25
  (1979)~(2) 175--178.

\bibitem{MR2519172}
G.~Wang and H.~Li, Color degree and alternating cycles in edge-colored graphs,
  Discrete Math. 309 (2009)~(13) 4349--4354.

\bibitem{Yeo1997}
A.~Yeo, A note on alternating cycles in edge-coloured graphs, J. Combin. Theory
  Ser. B 69 (1997)~(2) 222--225.

\end{thebibliography}
\end{document}